\documentclass[11pt,a4paper]{article}

\usepackage{lscape}
\usepackage{colortbl}
\usepackage{booktabs}
\usepackage[utf8]{inputenc}
\usepackage[T1]{fontenc}
\usepackage{authblk}
\usepackage[ backend=biber]{biblatex}
\usepackage{titlesec}
\titleformat{\paragraph}
{\normalfont\normalsize\bfseries}{\theparagraph}{1em}{}
\titlespacing*{\paragraph}
{0pt}{3.25ex plus 1ex minus .2ex}{1.5ex plus .2ex}
\usepackage{mathrsfs}
\usepackage[dvips]{epsfig}
\usepackage{graphicx}
\usepackage{amssymb, graphics, amsmath, amsthm, amsopn, amstext, amsfonts, algorithm, algorithmic}
\usepackage{subcaption}
\usepackage{fancyhdr}
\usepackage{lipsum} 
\usepackage{booktabs} 
\usepackage{tabularx}
\usepackage{caption}
\usepackage{adjustbox} 
\usepackage{rotating}
\usepackage{color}
\usepackage[left=3cm,right=2cm,top=3cm,bottom=2cm]{geometry}

\usepackage{multirow}
\usepackage{multicol}
\usepackage{cases}

\newcommand{\wt}{\widetilde}

\newcommand{\Argmin}{\mbox{Argmin}}

\newcommand{\Ical}{\mathcal I}
\newcommand{\Jcal}{\mathcal J}
\newcommand{\Lcal}{\mathcal L}

\newcommand{\mycut}[1]{{}}

\newtheorem{theorem}{Theorem}[section] 
\newtheorem{lemma}{Lemma}[section] 
\newtheorem{definition}{Definition}[section]

\makeatletter
\def\thanks#1{\protected@xdef\@thanks{\@thanks
        \protect\footnotetext{#1}}}
\makeatother

\begin{document}

\title{\textbf{
Cardinality Constrained Mean-Variance Portfolios: \\ A Penalty Decomposition Algorithm}}

\author[1]{Ahmad Mousavi$^*$\thanks{Emails: mousavi@american.edu and  gmichail@ucla.edu.
}} 
\affil{Department of Mathematics and Statistics, American University}
\author[2]{George Michailidis\thanks{}}
\affil[2]{Department of Statistics and Data Science, University of California, Los Angeles}

\date{}
\maketitle
\begin{abstract}
The cardinality-constrained mean-variance portfolio problem has garnered significant attention within contemporary finance due to its potential for achieving low risk while effectively managing risks and transaction costs. 
Instead of solving this problem directly, many existing methods rely on regularization and approximation techniques, which hinder investors' ability to precisely specify a portfolio's desired cardinality level. Moreover, these approaches typically include more hyper-parameters and increase the problem's dimensionality.{\color{black} \,
To address these challenges, we propose a customized penalty decomposition algorithm. We demonstrate that this algorithm not only does it converge to a local minimizer of the cardinality-constrained mean-variance portfolio problem, but is also computationally efficient. Our approach leverages a sequence of penalty subproblems, each tackled using Block Coordinate Descent (BCD). Specifically, we show that the steps within BCD yield closed-form solutions, allowing us to identify a saddle point of the penalty subproblems. Finally, by applying our penalty decomposition algorithm to real-world datasets, we highlight its efficiency and its superiority over state-of-the-art methods across several performance metrics.}
\end{abstract}

\section{Introduction}
Mean-variance portfolio optimization is of great importance in finance, enabling investors to balance risk and return through strategic diversification \autocite{kalayci2019comprehensive, mousavi2024statistical}. On the other hand, sparsity is an attractive property in mathematical optimization that has proved to be advantageous in various fields \autocite{mousavi2020survey}. For example, leveraging sparsity has proved useful in developing efficient algorithms, enabling for more concise and computationally economical representations, particularly in tasks involving large-scale data \autocite{moosaei2023sparse}. This property becomes particularly relevant when applied to mean-variance portfolios, where promoting sparsity significantly reduces transaction costs \autocite{lai2018short}.
Therefore, in this paper, we  study the following cardinality-constrained mean-variance (CCMV) portfolio problem \autocite{moreira2022alternating}:
\begin{equation}\label{pr: p-original} 
\quad \min_{x\in \mathbb{R}^n} \ f(x):= x^TAx - \tau \mu^Tx  \qquad \mbox{subject to}   \qquad  e^Tx=1, \quad  x\ge 0  \quad \text{and}\quad \|x\|_0\le k,
\end{equation}
where $\tau> 0$ balances the trade-off between minimizing risk and maximizing return, $\mu \in \mathbb R^n$ is the return vector, $A\succeq 0$ is the covariance matrix, and $k\in [n]$ is an upper bound for the number of desired assets. 

Since it is well-established that sparsity can be promoted using convex or nonconvex regularization techniques \autocite{mousavi2020survey}, many state-of-the-art studies resort to such approaches and solve approximations of (\ref{pr: p-original})  \autocite{kremer2020sparse, corsaro2019adaptive, zhao2021optimal}. However, this approach includes a regularization parameter that only implicitly controls the number of assets included in the resulting portfolio 
(that is, the number of nonzero elements in the obtained solution).
Further, this approach requires careful tuning of hyperparameters, which can be computationally expensive and challenging. Semidefinite relaxation methods have also been proposed \autocite{lee2020sparse}, but for medium or large-size problems are computationally expensive, since the number of decision variables in their formulation is squared. Moreover, at present, they lack rigorous theoretical justification that clarifies how to interpret the results of the corresponding semidefinite program based on the original problem. Note that (\ref{pr: p-original}) can be recast as a mixed integer programming similar to approaches employed in \autocite{bertsimas2009algorithm,gao2017dynamic}. Nevertheless, this thrust requires an exhaustive search, making it practically unsuitable for large-scale problems. Therefore, it is more advantageous to directly solve cardinality-constrained mean-variance portfolio problems, rather than their approximations or regularized variations, where the practitioner has the choice to explicitly specify the number of assets in the portfolio \autocite{shi2022cardinality}. 

To this point, for tackling (\ref{pr: p-original}) without the nonnegativity constraint, an ADMM method has been recently proposed \autocite{shi2022cardinality}. However, in addition to requiring penalty parameter tuning, this method incorporates the equality constraint into the objective function by means of another positive parameter that also must be selected carefully in practice. Further, the obtained solution of this method is not guaranteed to satisfy the first-order necessary optimality conditions. {\color{black}\autocite{moreira2022alternating} apply the generic methodology of Penalty Alternating Direction Method (PADM) to solve the cardinality-constrained mean-variance portfolio problem (\ref{pr: p-original}). CCMV-PADM involves solving a sequence of $\ell_1$ penalty subproblems using the BCD method. However, for each penalty parameter, this method necessitates solving a sequence of quadratic programs with a number of linear constraints equal to the problem size, which becomes computationally expensive for large-scale problems. Moreover, the theoretical foundation of CCMV-PADM is not robust, as the existence of a limit point is not guaranteed. Even assuming such a limit point exists, it is only a saddle point (see Definition \ref{def: saddle}).
}

To address these limitations, we present a customized penalty decomposition algorithm (CCMV-PD) for finding a local minimizer of the optimization problem (\ref{pr: p-original}). Our approach involves solving a series of $\ell_2$ penalty subproblems using the BCD method. Notably, each subproblem within the BCD step of our method has a closed-form solution, enhancing the computational efficiency of CCMV-PD compared to CCMV-PADM. Numerical evidence based on real data demonstrates the superior performance of the proposed algorithm over CCMV-PADM across several standard performance metrics and shows it is competitive with the commercial solver Mosek \autocite{cvx2013}. Additionally, we demonstrate that the CCMV-PD method outperforms both Mosek and CCMV-PADM in terms of CPU time.

\textbf{Notation}.
We define the set $[n]$ as $\{1,2,\dots,n\}$ for any natural number $n$. The set complement of $S$ is represented as $S^c$.
For a subset $S=\{i_1,i_2,\dots, i_{|S|}\}$, where $i_k$ are elements in $[n]$, and for any vector $x$ in $\mathbb{R}^n$, we define $x_S$ as the coordinate projection of $x$ with respect to the indices in set $S$. In other words, $(x_S)_i=x_i$ for $i$ belonging to $S$, and consequently, by slightly abusing the notation, an $S$-supported vector can be shown as $[x_S;0]$. We assess whether the matrix $A$ is positive semidefinite or positive definite using the symbols $A \succeq 0$ and $A\succ 0$, respectively. Also, $\lambda_{\min}(A)$ and $\lambda_{\max}(A)$ show the minimizer and maximum eigenvalue of $A$, respectively.

\section{A Penalty Decomposition Algorithm for CCMV and Its Convergence} \label{sec: algo}
Next, we introduce our customized penalty decomposition method (CCMV-PD) for solving the cardinality-constrained mean-variance problem (\ref{pr: p-original}) and establish that it obtains a local minimizer of this nonconvex program. 
\subsection{CCMV-PD Algorithm}
First, by introducing a new variable $y$, we equivalently reformulate (\ref{pr: p-original}) as follows:
\begin{equation} \label{pr: equivalent_p} 
\begin{aligned}
\min_{x,y\in \mathbb{R}^n} \quad 
x^TAx-\tau \mu^Tx \qquad
\textrm{subject to} \qquad 
 e^Tx=1, \quad y\ge 0, \quad
  \|y\|_0\le k, \quad \text{and} \quad     x-y=0.       \\
\end{aligned}
\end{equation}
This splitting effectively handles the nonlinearly coupled constraints of the original problem.
Suppose that
\begin{equation} \label{def: q_rho}
 q_{\rho}(x,y):= x^TAx-\tau \mu^Tx  +\rho\|x-y\|_2^2,
 \end{equation}
and 
\begin{equation*}
\mathcal X:=\{x\in \mathbb{R}^n \, | \, e^Tx=1\}, \qquad \text{ and } \qquad \mathcal Y:=\{y\in \mathbb{R}^n \, | \,y\ge 0 \quad \text{and} \quad  \|y\|_0\le k\}.
\end{equation*}
After penalizing the last constraint in (\ref{pr: equivalent_p}),  we consider a sequence of penalty subproblems as follows:
\begin{equation} \label{pr: pxy} 
\min_{x,y\in \mathbb R^n} \quad q_{\rho}(x,y)  \qquad
\textrm{subject to} \qquad 
 x\in \mathcal{X} \quad \text{and} \quad y\in \mathcal{Y}.
\end{equation}
To gain insights regarding the above formulation, it can be seen that by gradually increasing $\rho$ toward infinity, we can effectively address the optimization problem (\ref{pr: equivalent_p}). The introduction of the auxiliary variable $y$ is designed to simplify the subproblems encountered when implementing the subsequent algorithm for obtaining a saddle point in (\ref{pr: pxy}).

\begin{algorithm}[H]
\caption{BCD Algorithm for Solving (\ref{pr: pxy}) } 
\begin{algorithmic}[1]
\label{algo: BCD}
\STATE Input: Select arbitrary $y_{0}\in \mathcal Y$. 
\STATE Set $l=0$.
\REPEAT
\STATE 
$x_{l+1}=\Argmin_{x\in \mathcal X} \ q_{\rho}(x,y_l).$
\STATE 
$y_{l+1}= \Argmin_{y\in \mathcal Y} \ q_{\rho}(x_{l+1},y).$
\STATE $l \leftarrow l+1$ and go to step (3).
\UNTIL{stopping criterion (\ref{BCD-practical-stopping-criterion}) is met}. 
\end{algorithmic}
\end{algorithm}
Our stopping criterion for this BCD loop is the following:
\begin{equation} \label{BCD-practical-stopping-criterion}
\max
\left\{
\frac{\|x_l-x_{l-1}\|_\infty}{\max \left(\|x_l\|_\infty,1 \right)},
\frac{\|y_l-y_{l-1}\|_\infty}{\max \left(\|y_l\|_\infty,1 \right)}
\right\}
\le
\epsilon_I.
\end{equation}

We next study how to efficiently solve the restricted subproblems in Algorithm \ref{algo: BCD} below.
\newline
$\bullet
\ \Argmin_{x\in\mathcal X}\ q_{\rho}(x,y)$: This subproblem becomes the following convex quadratic optimization problem: 
\begin{equation} \label{pr: px}
\begin{aligned}
\min_{x\in \mathbb{R}^n} \quad & 
x^TAx-\tau \mu^Tx+\rho \|x-y\|^2
 \qquad \textrm{subject to} \qquad
 e^Tx=1.
\end{aligned}
\end{equation}
Since $A\succeq 0$ and $\rho>0$, we have $A+\rho I \succ 0$. Hence, this feasible problem has a solution and its KKT conditions are as follows:
\begin{equation*} \label{eqn: KKT _Px}
2 Ax_* -\tau \mu + 2\rho (x_*-y)+\beta e =0, \quad \mbox{ and } \quad e^Tx_*=1.
\end{equation*}
Thus,  $2(A+\rho I)x_*=\tau \mu +2\rho y -\beta e$ implies $x_*=0.5 (A+\rho I)^{-1}(  \tau \mu +2\rho y -\beta e)$. Since $1=e^Tx_*=0.5e^T(A+\rho I)^{-1}(  \tau \mu +2\rho y -\beta e)$, we get 
\begin{equation*}
\beta = \frac{-1+0.5e^T(A+\rho I)^{-1}(  \tau \mu +2\rho y)}{0.5e^T(A+\rho I)^{-1}e},
\end{equation*}
which in turn yields
\begin{equation} \label{eqn: x_*}
x_* = \frac{1}{2}(A+\rho I)^{-1}\left( 
\tau \mu +2 \rho y 
+
\frac{1-0.5e^T(A+\rho I)^{-1}(  \tau \mu +2\rho y)}{0.5e^T(A+\rho I)^{-1}e}e
\right)
\end{equation}
$\bullet
\ \Argmin_{y\in\mathcal Y}\ q_{\rho}(x,y)$: This subproblem is as follows:
\begin{equation} 
\label{pr: py} 
\min_{y\in \mathbb{R}^n} \ \|x-y\|_2^2 \qquad \mbox{subject to}  \qquad  y\ge 0 \quad \textrm{and} \quad \|y\|_0\le k.
\end{equation}
Recall that for $x\in \mathbb R$ and $k\in \mathbb N$, the sparsifying operator $H_k(x):=[x_\Jcal;0]$, where $\Jcal=\Jcal(x,k)$ is an index set corresponding to the $k$ largest components of $x$ in absolute value.
\begin{lemma} \label{lem: p_y solution}
The solution of problem (\ref{pr: py}) is  
\begin{equation} \label{eqn: y_*}
y_*= H_k(x^+)=H_k(\max(x,0)).
\end{equation}
\end{lemma}
{\color{black}
}
%
The CCMV starts with a positive penalty parameter and incrementally enlarges it until a convergence result is achieved. For a fixed $\rho$, Algorithm \ref{algo: BCD} handles the corresponding subproblem. Note that (\ref{pr: p-original}) is feasible, and we assume the availability of a feasible point denoted as $x^{\text{feas}}$.
To facilitate the presentation of this algorithm and its subsequent analysis, we introduce the following notation:
\begin{equation} \label{def: upsilon}
\Upsilon \ge \max\{f(x^{\text{feas}}), \min_{x\in \mathcal{X}} q_{\rho^{(0)}}(x,y^{(0)}_0) \} > 0 
\qquad \text{ and } \qquad 
X_\Upsilon := \{x \in \mathbb R^n\, | \, f(x)\le \Upsilon\}.
\end{equation}
Further, since $A\succeq 0,$ it is evident from a normalization argument that $X_\Upsilon$ is bounded and, therefore, compact.

\begin{algorithm}[H]
\caption{CCMV-PD Algorithm for Solving (\ref{pr: p-original}) }
\begin{algorithmic}[1]
\label{algo: CCMV-PD}
\STATE Inputs:  $\zeta>1, \rho^{(0)}>0$, and  $y^{(0)}_0\in \mathcal Y$.
\STATE Set $j=0$.
\REPEAT
 \STATE Set $l=0$.
 \REPEAT
\STATE Solve $x^{(j)}_{l+1} = \Argmin_{x\in \mathcal X} \ q_{\rho^{(j)}}(x,y^{(j)}_l).$
\STATE Solve
$y^{(j)}_{l+1} \in \Argmin_{y\in \mathcal Y} \ q_{\rho^{(j)}}(x^{(j)}_{l+1},y).$
 \STATE Set $l \leftarrow l+1$.
 \UNTIL{stopping criterion (\ref{BCD-practical-stopping-criterion}) is met.} \STATE Set
 $(x^{(j)},y^{(j)}):= (x^{(j)}_{l}, y^{(j)}_{l})$.
\STATE  Set $\rho^{(j+1)} = \zeta\cdot \rho^{(j)}$.
\STATE 
If $\min_{x\in \mathcal{X}} q_{\rho^{(j+1)}}(x,y^{(j)})> \Upsilon$, then $y_0^{(j+1)}=x^{\text{feas}}$.
Otherwise,
$y^{(j+1)}_0 = y^{(j)}$.

\STATE Set $j \leftarrow j+1$.
\UNTIL{stopping criterion (\ref{outer loop stopping criteria}) is met}. 
\end{algorithmic}
\end{algorithm}
The convergence of the outer  loop is met when:
\begin{equation} \label{outer loop stopping criteria}
\|x^{(j)}-y^{(j)}\|_{\infty}
\le
\epsilon_O.
\end{equation}

\subsection{Convergence Analysis of CCMV-PD Algorithm}
We begin by examining Algorithm \ref{algo: BCD} and subsequently establish the desirable properties of Algorithm \ref{algo: CCMV-PD}. 
We conduct an analysis of a sequence $\{(x_l,y_l)\}$ generated by  Algorithm~\ref{algo: BCD}, demonstrating that its accumulation point is a saddle point of (\ref{pr: pxy}), that is, we have:
$x_{*}\in \Argmin_{x\in \mathcal X} \ q_{\rho}(x,y_*),$ and $y_{*}\in\Argmin_{y\in \mathcal Y} \ q_{\rho}(x_*,y)$.
This shows the reason for choosing the BCD method for this nonconvex problem. 

{\color{black}
Within the convergence proof of Algorithm \ref{algo: CCMV-PD}, we need the boundedness of outer iterates $\{(x^{(j)},y^{(j)})\}$, which we establish below. This simply guarantees the boundedness of the BCD sequence for a fixed $\rho$, resulting in the existence of an accumulation point.
\begin{lemma}\label{lem: BCD boundedness}
Let $\{\rho^{(j)}\}$ be a positive monotone sequence going to infinity. Suppose $\{x^{(j)}\}$ and $\{y^{(j)}\}$ be two sequences such that $y^{(j)} \ge 0; \, \forall j$, and $\|x^{(j)}\|\ge \|y^{(j)}\|; \, \forall j$. Further, each $x^{(j)}$ is the unique minimizer of (4) for $\rho = \rho^{(j)}$ and $y=y^{(j)}$, that is,
\begin{equation*}
x^{(j)}= \frac{1}{2}(A+\rho^{(j)} I)^{-1}\left( 
\tau \mu +2 \rho^{(j)} y^{(j)} 
+
\frac{1-0.5e^T(A+\rho^{(j)} I)^{-1}(  \tau \mu +2\rho^{(j)} y^{(j)})}{0.5e^T(A+\rho^{(j)} I)^{-1}e}e
\right).
\end{equation*}
Then,  $\{x^{(j)}\}$ is bounded.
\end{lemma}
}
%
The next theorem proves that every accumulation point generated by the BCD is a saddle point of (\ref{pr: pxy}). Additionally, the sequence $\{q_\rho(x_l, y_l)\}$ is either strictly decreasing or two consecutive terms produce the same value, resulting in a saddle point. Essentially, Algorithm \ref{algo: BCD} produces a saddle point either in finite steps or in the limit.

\begin{theorem} \label{thm: BCD convergence}
Let $\rho > 0$. Consider the sequence $\{(x_l, y_l)\}$ generated by Algorithm~\ref{algo: BCD} to solve (\ref{pr: pxy}). If $(x_*, y_*)$ is an accumulation point of this sequence, then $(x_*, y_*)$ is a saddle point of the non-convex problem (\ref{pr: pxy}). Additionally, the sequence $\left\{q_{\rho}(x_{l}, y_{l})\right\}$ is non-increasing, and if $q_\rho(x_r, y_r) = q_\rho(x_{r+1}, y_{r+1})$ for some $r\in \mathbb N$, then $(x_r, y_r)$ is a saddle point of (\ref{pr: pxy}).

\end{theorem}

We now establish that Algorithm \ref{algo: CCMV-PD} obtains a local minimizer of (\ref{pr: p-original}). 
Assuming that $x^*$ is a local minimizer of (\ref{pr: p-original}), there exists an index set $\Lcal$ such that $|\Lcal|=k$ and $x^*_{\Lcal^c}=0$ such that $x^*$ is also a local minimizer of  the following problem:
\begin{equation}\label{pr: p-original-equivalent} 
\min_{x\in \mathbb{R}^n} \ x^TAx-\tau \mu^Tx \qquad \mbox{subject to}  \qquad  e^Tx=1, \quad  x\ge 0, \quad \text{and}\quad x_{\Lcal^c}=0.
\end{equation}
Robinson's constraint qualification condition for a local minimizer $x^*$ of (\ref{pr: p-original-equivalent}) is \autocite{mousavi2022penalty}:
\begin{equation} \label{eqn: robinson}
        \left\{  \begin{bmatrix} -d - v  
        \\ e^T d\\ d_{\Lcal^c} \end{bmatrix} \  \Big | \ \  d \in \mathbb R^{n}, \ \ v \in \mathbb R^n, \ \ v_i\le 0; \, \forall i\in \Lcal^c\right\} = \mathbb R^n \times \mathbb R \times \mathbb R^{|\Lcal^c|}. 
\end{equation}
It is easy to show that Robinson's conditions above always hold for an arbitrary $x^*$. Therefore, the KKT conditions for a local minimizer $x^*$ of (\ref{pr: p-original}) holds, which are the existence of $ (\lambda,\mu,w)\in \mathbb R^n \times \mathbb R\times \mathbb R^n $ and $\mathcal{ L}\subseteq [n]$ with $|\mathcal{L}|=k$ such that the following holds:
\begin{equation} \label{eqn: KKT _conditions for p}
2Ax-\tau \mu -\lambda +w+\beta e=0, \quad 0\le \lambda \perp x\ge 0, \quad e^Tx=1, \quad x_{{\mathcal{L}}^c}=0, \quad \mbox{and} \quad w_{\mathcal{L}}=0.
\end{equation}

Next, we analyze a sequence  $\{(x^{(j)}, y^{(j)})\}$ generated by Algorithm \ref{algo: CCMV-PD} and demonstrate that it finds a KKT  point of (\ref{pr: p-original}).

\begin{theorem} \label{thm: CCMV-PD convergence}
Suppose that  $\tau \in \mathbb R, k\in [n], A\succeq 0$, and $\rho^{(0)} \ge \lambda_{\max}(A)+1$.  
Let  $\{(x^{(j)}, y^{(j)})\}$  be a sequence generated by Algorithm~\ref{algo: CCMV-PD} for solving (\ref{pr: p-original}). Then, the following hold:
 \begin{itemize}
  \item [(i)] $\{(x^{(j)}, y^{(j)})\}$ has a convergent subsequence whose accumulation point  $(x^*, y^*)$ satisfies $x^*=y^*$. Further, there exists an index subset $\Lcal\subseteq [n]$ with $|\Lcal|=k$ such that $x^*_{\Lcal^c}=0$ and $x^*_{\Lcal}\ge 0$.
  \item [(ii)] Let $x^*$ and $\Lcal$ be defined above. Then, $x^*$ is a local minimizer of (\ref{pr: p-original}).
 \end{itemize}
\end{theorem}

{\color{black}
\subsection{A Comparison to CCMV-PADM}
To highlight the significance of the CCMV-PD algorithm in comparison with the CCMV-PADM one in \cite{moreira2022alternating} for solving (\ref{pr: p-original}), we begin with a preliminary discussion about the PADM methodology.

Consider the following optimization problem:
\begin{equation} \label{pr: two-variable-uv}
\min_{u,v} \ 
F(u,v) \qquad
\textrm{subject to} \qquad 
G(u,v)=0, \quad H(u,v)\ge 0, \quad u\in \mathcal{U},
\quad
\mbox{and}
\quad
v\in \mathcal{V}, 
\end{equation}
and let
$$
\Omega := \{(u,v) \in \mathcal{U} \times \mathcal{V} \, | \, G(u,v)=0 \, \, \text{and} \, \, H(u,v)\ge 0\}.
$$

\begin{definition}
A feasible point $(u^*, v^*) \in \Omega$ is called a saddle point of (\ref{pr: two-variable-uv}) if
\begin{equation} \label{def: saddle}
F(u^*, v^*) \le F(u, v^*), \quad \forall (u, v^*) \in \Omega \quad \text{and} \quad F(u^*, v^*) \le F(u^*, v), \quad \forall (u^*, v) \in \Omega,
\end{equation}
which is equivalent to $u^* \in \Argmin_{u \in \mathcal{U}} F(u, v^*)$ and $v^* \in \Argmin_{v \in \mathcal{V}} F(u^*, v),$ respectively.
\end{definition}

The feasible set $\Omega$ is often not decomposable into disjoint sets $\mathcal{U}$ and $\mathcal{V}$ due to the presence of coupling constraints $G(u,v)=0$ and $H(u,v) \ge 0$. This interdependence results in poor convergence of the standard BCD method in practice. To address this issue, PADM introduces the following $\ell_1$ penalty function:

\begin{equation} \label{pr: phi1}
\phi_1(u,v; \rho, \theta) = F(u,v) + \sum_{i}\rho_i |G_i(u,v)| + \sum_{i} \theta_i \max(-H_i(u,v), 0),\end{equation}
with $\rho_i \ge 0$ and $\theta_i \ge 0; \forall \ i,$ being the penalty parameters. The PADM starts with initial positive values for penalty parameters and solves a sequence of penalty subproblems as follows:
\begin{equation*}
\min_{u,v}
\, \phi_1(u,v; \rho, \theta) \qquad
\text{subject to} \qquad u \in \mathcal{U},
\quad \mbox{and} \quad
v \in \mathcal{V},
\end{equation*}
where the BCD method is utilized to find a saddle point of the above. The PADM stops once an obtained saddle point belongs to $\Omega$, otherwise, the penalty parameters are increased and the process is repeated. To sum it up, the PADM generates a sequence of saddle points (each corresponding to a fixed penalty parameters pair using the BCD method) until one of them satisfies the coupling constraints $G(u,v)=0$ and $H(u,v)\ge 0$.
The main convergence result of this method is reported below.

\begin{theorem} \label{thm: PADM}
Suppose that functions $F, G,$ and $H$ are continuous and $\mathcal{U}$ and $\mathcal{V}$ are nonempty and compact. Let $\rho^{(j)} \nearrow \infty$ and $\theta^{(j)} \nearrow \infty$ elementwise. Let $\{(u^{(j)},v^{(j)})\}$ be a sequence of saddle points of (\ref{pr: phi1}) corresponding to $\rho^{(j)}$ and $\theta^{(j)}$ and $(u^{(j)},v^{(j)})\to (u^*,v^*)$. Let $(u^*,v^*)$ belong to $\Omega$. Then, $(u^*,v^*)$ is saddle point of (\ref{pr: two-variable-uv}). Further, if $F$ is continuously differentiable, then  $(u^*,v^*)$ is a KKT point of (\ref{pr: two-variable-uv}). 
\end{theorem}

To apply this generic methodology on (\ref{pr: p-original}), \cite{moreira2022alternating} add a copy constraint to it as follows:
\begin{equation}\label{pr: p-original-morrrrr} 
\quad \min_{x,y\in \mathbb{R}^n} \ F(x,y):= x^TAx - \tau \mu^Tx  \qquad \mbox{subject to}   \qquad  e^Tx=1 \quad  x\ge 0, \quad x=y, \quad \text{and}\quad \|y\|_0\le k.
\end{equation}
Then, they consider the following penalty subproblem:
\begin{equation*} 
\begin{aligned}
\min_{x,y\in \mathbb{R}^n} \ & 
f(x,y):= 
x^TAx - \tau \mu^Tx+\rho \|x-y\|_1 \\
\mbox{subject to}  \ \ & x\in \tilde{\mathcal{U}}:=\{x \, | \, e^Tx=1 \ \ \mbox{and} \, \, x\ge 0\}  \quad \mbox{and} \quad
w\in \tilde{\mathcal{V}}:=\{w \, | \, e^Ty=1 \, \, \mbox{and} \, \, \|y\|_0\le k\}.
\end{aligned}
\end{equation*}
This means that those subproblems that emerge in the BCD step of the penalty alternating direction method (CCMV-PADM) given in \cite{moreira2022alternating} are as follows:
\begin{equation}  \label{pr: px-morierra}
\begin{aligned}
\min_{x\in \mathbb{R}^n} \quad 
x^TAx-\tau \mu^Tx+\rho\|x-y\|_1\qquad
\textrm{subject to} \qquad 
 e^Tx=1 \quad \text{and}  \quad x\ge 0,
\end{aligned}
\end{equation}
and
\begin{equation} \label{pr: py-morierra}
\begin{aligned}
\min_{y\in \mathbb{R}^n} \quad 
\|x-y\|_1\qquad
\textrm{subject to} \qquad 
 e^Ty=1 \quad \text{and}  \quad \|y\|_0\le k.
\end{aligned}
\end{equation}
Several issues with the CCMV-PADM algorithm are discussed next: 
\begin{itemize}
    \item[(i)]
Although \cite{moreira2022alternating} demonstrates that (\ref{pr: py-morierra}) has a closed-form solution, (\ref{pr: px-morierra}) is a nondifferentiable convex program, equivalent to a quadratic program with $n$ constraints, where $n$ is the problem size. Consequently, solving such a subproblem iteratively increases computational cost. In contrast, our proposed algorithm is computationally more efficient as both of our subproblems have closed-form solutions. Further, (\ref{pr: px-morierra}) lacks a unique solution, which may lead to weak convergence in practice.
\item[(ii)]
The sequence generated by CCMV-PADM may not be bounded, specifically because the set $\tilde{\mathcal{V}}$ in this method is unbounded. Therefore, the existence of an accumulation point is uncertain. However, \cite{moreira2022alternating} did not encounter practical issues and observed convergence in their numerical results. Nevertheless, a theoretical proof is still necessary. In contrast, Lemma \ref{lem: BCD boundedness} ensures the existence of an accumulation point for any sequence generated by CCMV-PD.

\item[(iii)]
Even if the existence of an accumulation point is assumed in CCMV-PADM, note that $f(x,y)$ is not continuously differentiable. Consequently, a limit point in this method is only guaranteed to be a saddle point of (\ref{pr: p-original-morrrrr}) or equivalently (\ref{pr: p-original}). In contrast, Theorem \ref{thm: CCMV-PD convergence} guarantees convergence to a local minimizer of (\ref{pr: p-original}).

\end{itemize}

}

\section{Numerical Results} \label{sec: numerical}
In this section, we first compare the performance of CCVM-PD (\ref{algo: CCMV-PD}) for solving (\ref{pr: p-original}) with CCMV-PADM proposed in \autocite{moreira2022alternating} and the commercial solver Mosek available in CVX \autocite{grant2014cvx} in terms of in-sample and out-of-sample returns, risks, and Sharpe ratios over the S\&P index\footnote{https://finance.yahoo.com}. We consider the period $2016-2020$ and randomly select $133$ stocks. For the out-of-sample performance, we apply the rolling-horizon procedure  \autocite{dai2018generalized}.
The formula for calculating the in-sample Sharpe ratio is ${SR} ={\mu}/{\sigma},$
where $ \mu $ is expected portfolio return and  $ \sigma $ is the standard deviation of  portfolio
and the out-of-sample Sharpe ratio ($\widehat{SR} $) is as follows:
\begin{align*}
	\widehat{SR} =\dfrac{\widehat{\mu}}{\widehat{\sigma}}, \qquad \text{ with } \qquad
	\widehat{\mu}=\dfrac{1}{T-\nu} \sum_{t=\nu}^{T-1} ~(x _{t}^T~r_{t+1})
 \qquad \text{ and } \qquad 
	\widehat{\sigma}^2 =\dfrac{1}{T-\nu -1} \sum_{t=\nu}^{T-1} ~(x _{t}^T~r_{t+1} - \widehat{\mu})^2,
\end{align*}
where  $  x _{t} $ denotes the stock weight at time $ t$, and   $ T $ is the total number of returns in the dataset, $ r_{t+1} $ is the stock return, and $ \nu $ is the length of the estimation time window.
By considering the monthly return of stocks from $2016$ to $ 2019 $ and an estimation window of $ \nu= 48 $ data, for $ 4 $ years, we use the 2020 data for the out-of-sample.
The results are summarized in Table \ref{table1} for when $ n=133, \tau =0.5,  \rho_0 = 0.1, \zeta=10 $ and different  values of $k$ are considered. 
We compare the performance of  Algorithm \ref{algo: CCMV-PD} in terms of returns, risks, and Sharpe ratios for the model (\ref{pr: p-original}) with  Mosek and CCMV-PADM for the dataset of  S\&P index from  2016 until 2019  when $ T=60 $. 
The BCD loop is stopped when  $\epsilon_I=10^{-4}$ in (\ref{BCD-practical-stopping-criterion}). Algorithm \ref{algo: CCMV-PD} also terminates  if $\epsilon_o=10^{-4}$ in (\ref{outer loop stopping criteria}). 
The results are summarized in Table \ref{table1}. In this table, the gaps in the last three columns are $ {|g-\widehat{g}|}/({|\widehat{g}|+1}),$ where $ g$ and $ \widehat{g}$  are the return, risk, Sharpe ratio Mosek. 

\begin{figure}[h] 
\centering
\includegraphics[height=11cm, width=.7\linewidth]{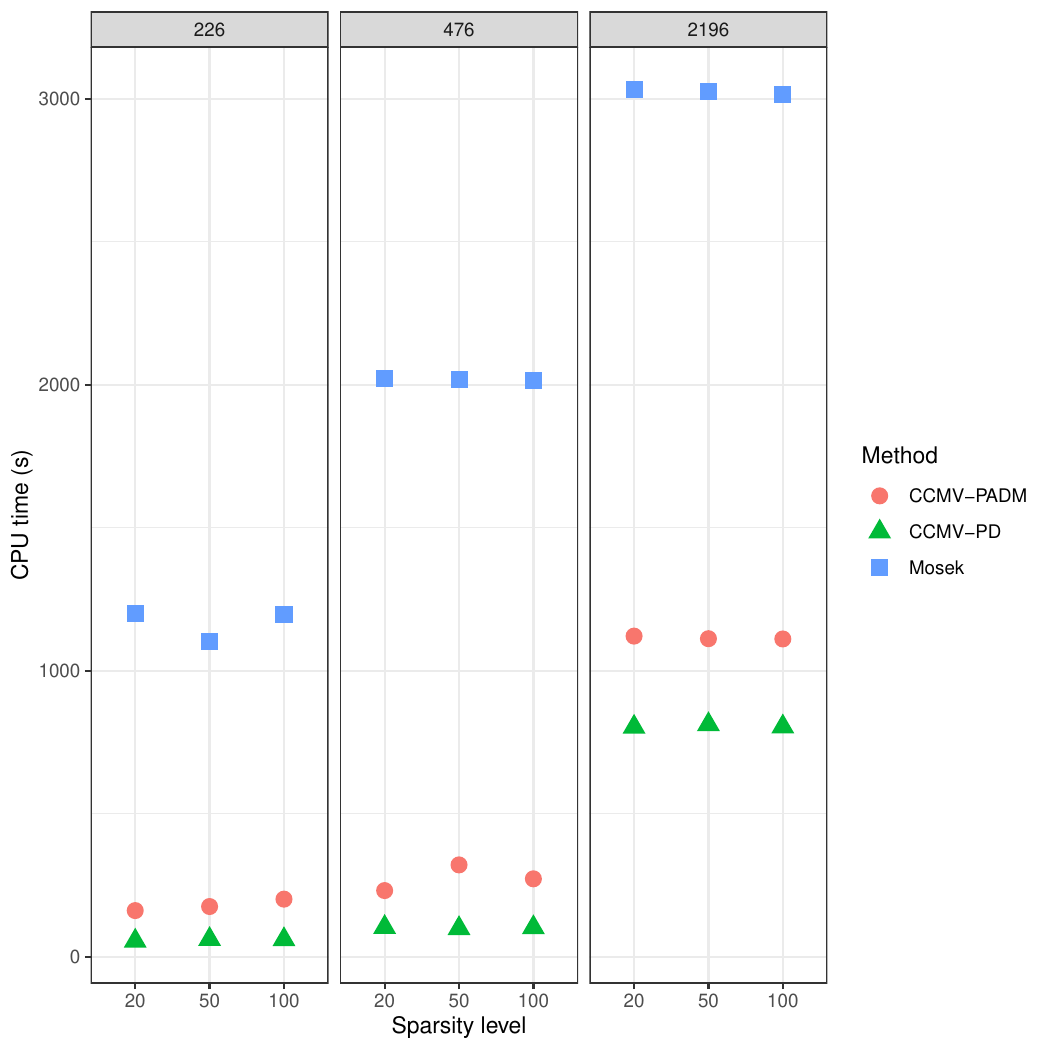}
\caption{Comparison of log of CPU times in Table \ref{table2} against sparsity level for $n=226, 476,$ and $2196$.}
\label{fig2}
\end{figure}

Our results reveal that  CCMV-PD and CCMV-PADM overall have close returns, risks, and Sharpe ratios, however, our proposed method does significantly better in a time comparison. Further, for small sparsity levels, this closeness of Mosek and CCMV-PD is significant. These results can be further confirmed in the graphs of Figure \ref{ek}.
As one can see, from the results of out-of-sample in Table \ref{table1}, the return, risk, and Sharpe ratio of all used methods in solving (\ref{pr: p-original}) are close to each other, which can be seen in the small gaps of these models.
Secondly, we compare the performance of Algorithm \ref{algo: CCMV-PD} in terms of CPU time for the model (\ref{pr: p-original}) with Mosek and CCMV-PADM on the dataset of  MIBTEL,  S\&P, and  NASDAQ indexes for different values of $k$. These results are reported in Figure \ref{fig2} and Table \ref{table2} and show that CCMV-PD solves (\ref{pr: p-original}) significantly faster and has smaller gaps. Due to the lack of space in this table, we shortened the return gap to re gap, risk gap to ri gap,  and Sharpe ratio to Sr. We conclude that the results of CCMV-PD are competitive to the results of the solving model (\ref{pr: p-original}) by Mosek in CVX. Further, our method shows a better performance than CCMV-PADM overall. 

\begin{landscape}
\begin{figure} [h]
\centering
\begin{minipage}{ .45\textwidth}
  \centering  \includegraphics[width=1\linewidth]{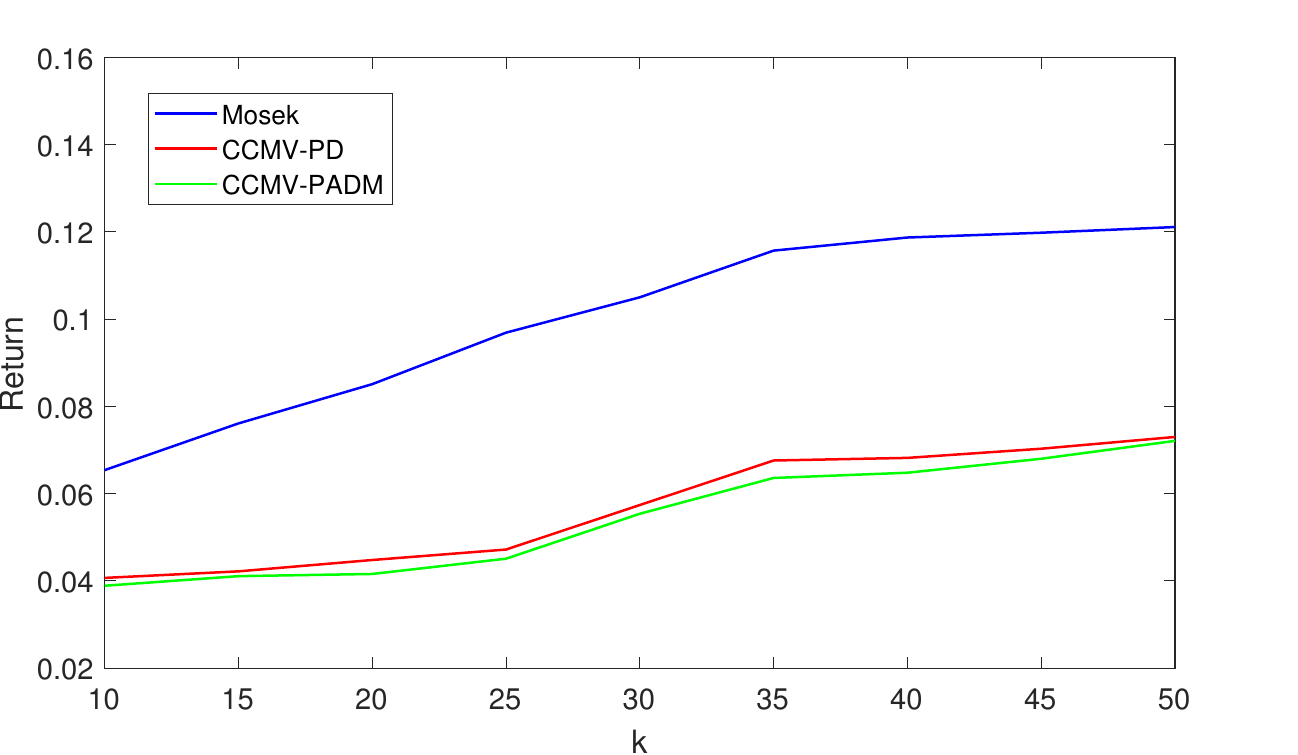}
\end{minipage}%
\begin{minipage}{ .45\textwidth}
  \centering
  \includegraphics[width=1\linewidth]{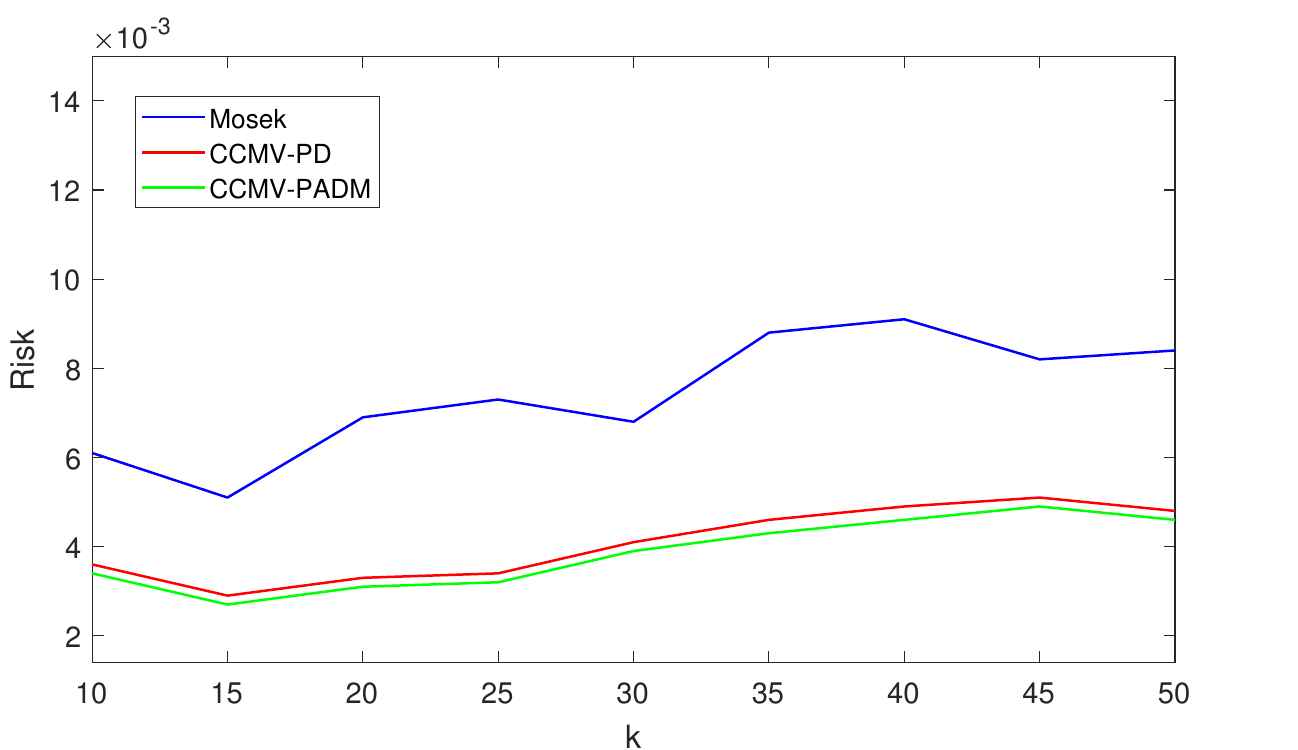}
\end{minipage}
\begin{minipage}{ .45\textwidth}
  \centering
  \includegraphics[width=1\linewidth]{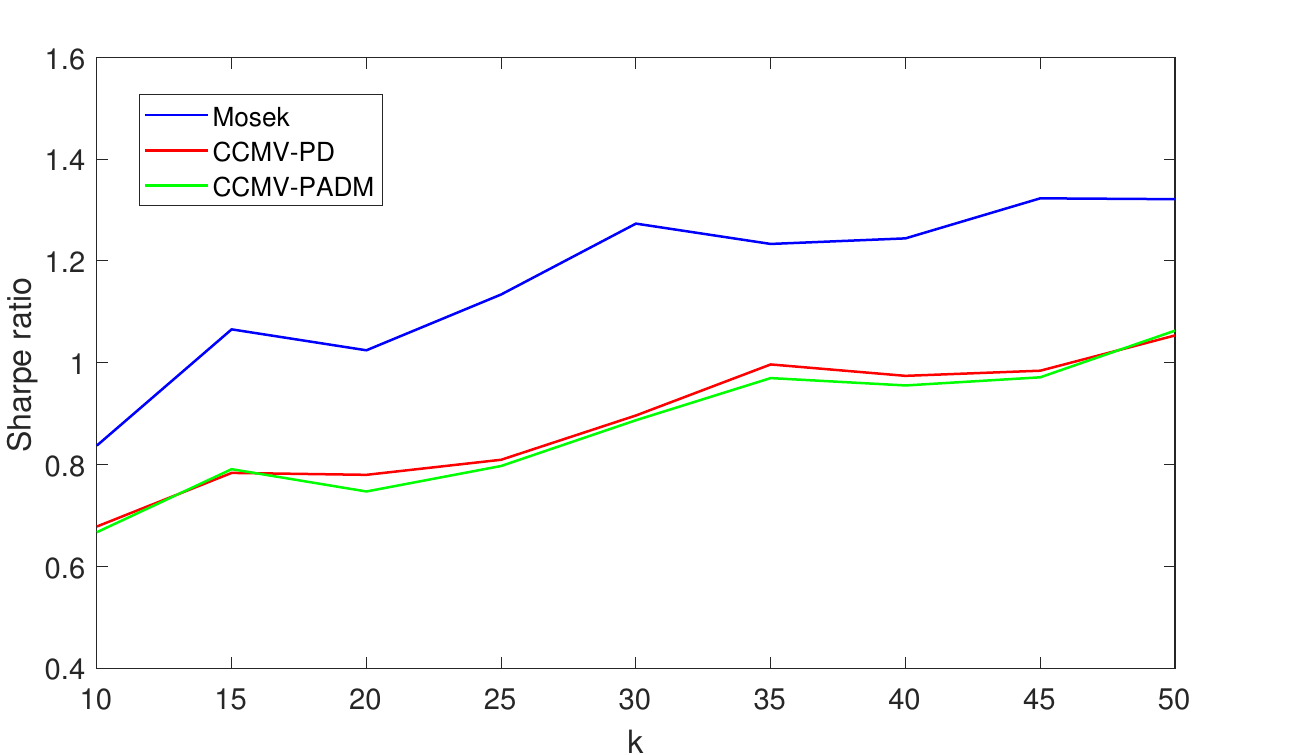}
\end{minipage}
\caption{Gap returns, risks, and Sharpe ratios for Mosek, CCPV-PD, and CCMV-PADM in solving (\ref{pr: p-original})  against sparsity level $k$.}
\label{ek}
\end{figure}
\bigskip \bigskip \bigskip 
\bigskip \bigskip 
\begin{table}[h] 
\begin{tabular}{cccccccccccccccccccccc}
\toprule
\multicolumn{1}{c}{} & 
\multicolumn{3}{c}{\textbf{Mosek}} & 
\multicolumn{3}{c}{\textbf{CCMV-PD}} &
\multicolumn{1}{c}{} &
\multicolumn{1}{c}{} &
\multicolumn{1}{c}{} &
\multicolumn{3}{c}{\textbf{CCMV-PADM}} &
\multicolumn{1}{c}{} &
\multicolumn{1}{c}{} &
\multicolumn{1}{c}{} 
\\
\cmidrule(rl){2-4} \cmidrule(rl){5-7} \cmidrule(rl){11-13}
 $k$ & return & risk & Sr &return & risk & Sr & re gap & ri gap & Sr gap  &return & risk & Sr  & re gap & ri gap & Sr gap \\
\midrule
20 & 0.0840  &  0.0051& 1.1765 & 0.0663 &  0.0038 & 1.0755 & 0.0163 &  0.0013 & 0.0464& 0.0654 &  0.0038 & 1.0609 & 0.0172 &   0.0013& 0.0531 \\
30 & 0.0808 & 0.0078 & 0.9151 & 0.0477 &0.0020 & 0.6194 & 0.0306 &  0.0058 & 0.1544 & 0.0372&  0.0038 &  0.4412 & 0.2500 &0.0020 & 0.2475 \\
40 & 0.0838  &  0.0053&  1.1512 &  0.0441 &  0.0024 & 0.9002 & 0.0366 & 0.0028 & 0.1167 & 0.0428 &  0.0029 & 0.6035 & 0.0378 & 0.0024&  0.1657  \\
50 &  0.1087 & 0.0088  & 1.1589 &  0.0456 &  0.0027  &  0.8776 & 0.0569 & 0.0061 & 0.1303 & 0.0448 &  0.0031 & 0.8046 & 0.0576 &  0.0057 & 0.1641  \\
20 & 0.1615 &  0.0082  &   1.7835 & 0.0653 &  0.0035 & 1.1038 & 0.0828 &  0.0048  &  0.2442 & 0.0651& 0.0035 &  1.1004 & 0.0830 &  0.0047 & 0.2454  \\
30 &	0.1609  &  0.0068  & 1.9512  & 0.0651 &  0.0031 &  1.1694 &  0.0825 &  0.0037  & 0.2649 & 0.0646 &  0.0031 & 1.1603 & 0.0829& 0.0037 & 0.2680  \\
40 & 0.1831  &  0.0087  &  1.9630 & 0.0669 &  0.0037 & 1.0998 & 0.0982  &  0.0051  &  0.2913 & 0.0651 & 0.0039 & 1.0424 & 0.0997 & 0.0048 & 0.3107  \\
50 & 0.1836  &  0.0089  &  1.9462  & 0.0678 & 0.0041  & 1.0589 & 0.0978  &  0.0049  &  0.3012 & 0.0666 & 0.0040  & 1.0530 & 0.0989  & 0.0049 & 0.3032 \\
\bottomrule
\end{tabular}
\caption{Comparison of returns, risks, and Sharpe ratios of Mosek, CCMV-PD, and CCMV-PADM for different sparsity levels when $n=133$} \label{table1}
\end{table}
\end{landscape}

\begin{landscape}
\begin{table}[h] 
\centering
\begin{tabular}{cccccc}
\toprule
\multicolumn{1}{c}{} & 
\multicolumn{1}{c}{} & 
\multicolumn{4}{c}{\textbf{Mosek}} 
\\
\cmidrule(rl){3-6}
 & $k$ & return & risk & Sr & time  \\
\midrule
& 20 & 0.0341 & 0.0028& 0.6364 & 1.201e+03  \\
$n=226$&  50 & 0.0368 & 0.0035& 0.6221 & 1.101e+03 
\\
&  100 & 0.0432 & 0.0051& 0.6049& 1.1982e+03  \\
& 20 & 0.0328  & 0.0030  & 0.5965  &  $ > $2.021e+03 
\\
$n=476$&  50  & 0.0550 &  0.0045 & 0.8213 &  $ > $2.018e+03  
\\
&  100  & 0.0583 & 0.0052 &  0.8085 &  $ > $2.014e+03 
\\
&  20  & 0.0615 &  0.0065 & 0.7631  &  $ > $3.031e+03  
\\
$n=2196$&  50  & 0.0623 &  0.0066 & 0.7669  &  $ > $3.024e+03 
\\
&  100  & 0.0645 &  0.0069 &  0.6515   &  $ > $3.015e+03 
\\
\bottomrule
\bigskip \bigskip
\end{tabular}
\begin{tabular}{cccccccccccccccc}
\toprule
\multicolumn{4}{c}{\textbf{CCMV-PD}} &
\multicolumn{1}{c}{} &
\multicolumn{1}{c}{} &
\multicolumn{1}{c}{} &
\multicolumn{4}{c}{\textbf{CCMV-PADM}} &
\multicolumn{1}{c}{} &
\multicolumn{1}{c}{} &
\multicolumn{1}{c}{} 
\\
\cmidrule(rl){1-4} \cmidrule(rl){8-11} 
return & risk & Sr & time &  re gap & ri gap & Sr gap  &return & risk & Sr  & time & re gap & ri gap & Sr gap \\
\midrule
 0.0098  &  4.9697e-04&0.4434  & 55.1256 & 0.0234 & 0.0023 & 0.1179 & 0.0091 & 4.3275e-04 & 0.4374 & 161.324 & 0.0242 & 0.0024 & 0.1216 \\
  0.0125  &  4.7697e-03  & 0.1810 & 61.0236 & 0.0177 &  0.0014 &  0.2719 &  0.0118 & 4.2994e-03 & 0.1800 & 175.454 & 0.0224 & 7.9661e-04 &  0.2725  
\\
0.0187  & 1.6697e-03 & 0.4576  & 60.3265 & 0.0234  &  0.0034  & 0.0918 &0.0181 &2.2881e-04 &0.4562 & 201.345 &0.0241 &  0.0048 & 0.0927  \\
0.0208  & 4.5780e-03& 0.3074 & 103.1256 & 0.0116 & 0.0016 & 0.1811  & 0.0198 & 4.1776e-03 & 0.3063 & 231.304&  0.0126 & 0.0012  & 0.1818 
\\
 0.0215 & 0.0065 & 0.2667 & 99.1325 & 0.0318 & 0.0020&  0.3045 &0.0201 & 0.0061 & 0.2574 & 321.104 & 0.0331 & 0.0016 & 0.3096
\\
0.0288 & 6.6639e-03 &  0.3528  & 102.1595 & 0.0279 & 0.0015 & 0.2520 &0.0281 & 6.3839e-03 & 0.3517  & 272.326 &  0.0285 & 0.0013 & 0.2526 
\\
 0.0266  & 0.0098 & 0.2687 & 803.1132 &  0.0329 & 0.0033& 0.2804 & 0.0256  & 0.0096 & 0.2613 & 1.121e+03 & 0.0338 &0.0031 & 0.2846
\\
0.0267& 9.8785e-03 & 0.2686 & 812.1335 & 0.0468 &0.0033 & 0.2820 & 0.0105 & 9.5669e-03 & 0.2678 & 1.112e+03& 0.0478 &  0.0029 &0.2825
\\
0.0123  & 9.8955e-03 &  0.4762 & 804.1635 &  0.0490 & 0.0092  &  0.1061 & 0.0118 & 9.2697e-04 &  0.4522 & 1.111e+03& 0.0495& 0.0092 &  0.1207
\\
\bottomrule
\end{tabular} \caption{Comparison of returns, risks, Sharpe ratios, and CPU times of Mosek, CCMV-PD, and CCMV-PADM in solving (\ref{pr: p-original})  for different numbers of stocks: $226$ stocks of MIBTEL, $476$ stocks of S\&P, and $2196$ stocks of  NASDAQ  indexes.} \label{table2}
\end{table}
\end{landscape}

{\color{black}
\section{Conclusion}
This paper addressed the cardinality-constrained mean-variance portfolio problem (\ref{pr: p-original}), a topic of significant interest in contemporary finance due to its potential for achieving low risk while effectively managing transaction costs and other risks. The proposed algorithm directly solves this cardinality problem, rather than its regularized versions or relaxations, thus enabling investors to precisely specify the desired cardinality level of their portfolios.

We introduced CCMV-PD, a customized penalty decomposition algorithm that fully exploits the structure within (\ref{pr: p-original}). CCMV-PD not only guarantees convergence to a local minimizer of the cardinality-constrained mean-variance portfolio problem, but also offers computational efficiency. Specifically, CCMV-PD solves a sequence of $\ell_2$ penalty subproblems, each approximately solved via the BCD method. We demonstrated that each step within the BCD method finds closed-form solutions and identifies a saddle point of the penalty subproblems. Ultimately, we established that these saddle points have an accumulation point that is indeed a local minimizer of (\ref{pr: p-original}). 

Further, the paper discusses the theoretical significance and strengths of CCMV-PD vis-a-vis the state-of-the-art algorithm CCMV-PADM presented in \cite{moreira2022alternating}, which employs a similar methodology by solving a sequence of $\ell_1$ penalty subproblems. Further, through extensive numerical experiments on real-world datasets, we showcased the practical efficacy of our algorithm. The results highlighted the superiority of CCMV-PD over state-of-the-art methods across several performance metrics, such as in-sample and out-of-sample returns, risks, and Sharpe ratios. This demonstrated the comprehensive advantages of CCMV-PD.

In summary, our proposed CCMV-PD algorithm provides a robust solution to the cardinality-constrained mean-variance portfolio problem, improving upon existing methods in both theoretical guarantees and practical performance. Future work may explore extending this framework to other types of portfolio optimization problems and further enhancing its computational efficiency.
}

\printbibliography

\section{Appendix} \label{sec: appendix}
Technical proofs are given below.

\subsection{Proof of Lemma \ref{lem: p_y solution}}

\begin{proof}
Note that for any $y$ satisfying $\| y \|_0 \le k$, one can write $y=(y_\Ical, y_{\Ical^c})$ such that $y_\Ical =0$ for some index set $\Ical \subseteq \{1, \ldots, n \}$ with $|\Ical|=n-k$. Hence, for any index set $\Ical$ with $|\Ical| =n-k$, $(P_y)$ can be written as 
$$\min_{ y \in \mathbb R^n} (\|  x_\Ical-y_\Ical \|^2_2+ \|  x_{\Ical^c}-y_{\Ical^c} \|^2_2) \quad \text{subject to } \quad y\ge 0 \quad \text{and} \quad y_{\Ical}=0,$$
which is equivalent to $$\min_{z \in \mathbb R^{|\Ical^c|} } \| z- x_{\Ical^c} \|^2_2 \quad \text{subject to} \quad z\ge 0.$$
Clearly, constraint qualification holds, and its KKT condition for a local minimizer $z_*$ is: $ x_{\Ical^c} =z_* -\lambda $ and $z_*\ge 0$ for some $\lambda \ge 0$ with $z_*^T\lambda =0$. 
Suppose that $x_{\Ical^c}=x^+_{\Ical^c}-x^-_{\Ical^c}$ such that $x^+_{\Ical^c}=\max(x_{\Ical^c},0)$ and $x^-_{\Ical^c}=\max(-x_{\Ical^c},0)$, then $x^{+}_{\Ical^c}\perp x^-_{\Ical^c}$. It is easy to see that $z_*=x^+_{\Ical^c}$ and $\lambda = x^-_{\Ical^c}$ satisfy the KKT conditions. Hence, it is easy to show that $z_* = x^+_{\Ical^c}$ for any index $\Ical$ specified above. Finally, for any index $\Ical$ specified above, the optimal value is given by
$ \|x_{\Ical^c} -  x^+_{\Ical^c}  \|^2_2 + \|x_\Ical\|_2^2 = \|x^-_{\Ical^c} \|^2_2 + \|x_\Ical\|^2_2=\|x_{\Ical^c}\|_2^2 -  \|x^+_{\Ical^c}  \|^2_2 + \|x_\Ical\|^2_2= \|x\|_2^2- \|x^+_{\Ical^c}  \|^2_2$. Consequently, the minimal value of $(P_y)$ is achieved when $\|x^+_{\Ical^c}\|_2$ is maximal or equivalently when $\Ical^c = \Jcal(x, k)=\Jcal$. Therefore, a minimizer $y^*$ satisfies $y_\Jcal = H_k(x^+)=H_k(\max(x,0))$.
\end{proof}

\subsection{Proof of Lemma \ref{lem: BCD boundedness}}

{\color{black}
\begin{proof} 
Suppose not.  Then, there is a subsequence $\{x^{(j_k)}\}$ of $\{x^{(j)}\}$ such that $\|x^{(j_k)}\|\rightarrow\infty$ as $j\rightarrow\infty$. Since $\left\|\frac{x^{(j_k)}}{\|x^{(j_k)}\|}\right\|=1,$$\left\|\frac{y^{(j_k)}}{\|x^{(j_k)}\|}\right\|\le1,$ the sequence $\left\{\left(\frac{x^{(j_k)}}{\|x^{(j_k)}\|},\frac{y^{(j_k)}}{\|x^{(j_k)}\|}\right)\right\}$ lies inside a closed ball of radius $\sqrt{2}$ of $\mathbb{R}^{2n};$ a compact set. Therefore, there is a subsequence $\left\{\left(\frac{x^{(j_l)}}{\|x^{(j_l)}\|},\frac{y^{(j_l)}}{\|x^{(j_l)}\|}\right)\right\}$ of  
$\left\{\left(\frac{x^{(j_k)}}{\|x^{(j_k)}\|},\frac{y^{(j_k)}}{\|x^{(j_k)}\|}\right)\right\}$ that converges, that is, 
\begin{align*}
\lim_{j\rightarrow\infty}\frac{x^{(j_l)}}{\|x^{(j_l)}\|}=\hat{x}\qquad \text{ and } \qquad
\lim_{j\rightarrow\infty}\frac{y^{(j_l)}}{\|x^{(j_l)}\|}=\hat{y},
\end{align*}
where $1=\|\hat{x}\|\ge\|\hat{y}\|,$ and $\hat{y}\ge0.$
For any $\rho>0$, we have:
$
(A+\rho I)^{-1}=\rho^{-1}({\rho}^{-1}A+I)^{-1},
$
we get
$
\lim_{\rho^{(j)}\rightarrow\infty}(\rho^{-1}A+I)^{-1}=I.
$
Hence, using (7), we get
\begin{align*}
x^{(j_l)}=\frac{1}{2}({\rho^{(j_l)}}^{-1}A+I)^{-1}\big({\rho^{(j_l)}}^{-1}\tau\mu +2y^{(j_l)}+
\frac{{\rho^{(j_l)}}^{-1}-0.5e^T({\rho^{(j_l)}}^{-1}A+I)^{-1}({\rho^{(j_l)}}^{-1}\tau\mu +2y^{(j_l)})}{0.5e^T({\rho^{(j_l)}}^{-1}A+I)^{-1}e}e
\big).
\end{align*}
Dividing this identity through by $\|x^{(j_l)}\|$ and taking limits as $j\rightarrow\infty$, gives
$$
\hat{x}=\hat{y}-\frac{e^T\hat{y}}{n}e.
$$
Taking norms from both sides implies
\begin{align}
1&=\|\hat{x}\|^2\notag \\
&=\|\hat{y}\|^2-\frac{(e^T\hat{y})^2}{n}\notag\\
&\le1-\frac{(e^T\hat{y})^2}{n}\notag,
\end{align}
which is impossible, because if $\hat{y}\ne0$ then $e^T\hat{y}>0$ and therefore $1-\frac{(e^T\hat{y})^2}{n}<1,$ and if $\hat{y}=0$ it gives $1=0.$ Therefore $\{x^{(j)}\}$ must be bounded.
\end{proof}
}


\subsection{Proof of Theorem \ref{thm: BCD convergence}}
\begin{proof}

By observing definitions of $x_{l+1}$ and $y_{l+1}$ in steps 4 and 5 of Algorithm \ref{algo: BCD}, we get
\begin{align} \label{eqn: q_xy-inequality}
  q_{\rho}(x_{l+1},y_{l+1}) &\le  q_{\rho}(x_{l+1},y), \quad \forall y\in \mathcal Y, \notag \\
  q_{\rho}(x_{l+1},y_{l}) &\le  q_{\rho}(x,y_{l}), \quad \forall x\in \mathcal X.
\end{align}
This simply leads to the following:
\begin{equation} \label{eqn: non-increasing_q}
\begin{aligned}
c^* &\le  q_{\rho}(x_{l+1},y_{l+1}) \le q_{\rho}(x_{l+1},y_{l})\le q_{\rho}(x_{l},y_{l}), \quad \forall l \in \mathbb N,
\end{aligned}
\end{equation}
where $c^*=\min_x x^TAx-\tau \mu^Tx,$ which is finite because $A\succ 0$. 
Thus, $q_{\rho}(x_{l},y_{l})$ is a bounded below and non-increasing sequence; implying that $q_{\rho}(x_{l},y_{l})$ is convergent. 
From the other side, since $(x_*,y_*)$ is an accumulation point of $\{(x_l,y_l)\}$, there exists a subsequence $L$ such that $\lim_{l\in L\to \infty} (x_l,y_l) = (x_*,y_*)$. The continuity of $q_{\rho}(x_{l},y_{l})$ yields 
\begin{equation*} 
\begin{aligned}
\lim_{l\to \infty}q_{\rho}(x_{l+1},y_{l+1}) &= \lim_{l\to \infty} q_{\rho}(x_{l+1},y_{l}) \\
&=  \lim_{l\to \infty} q_\rho(x_l,y_l) \\
&= \lim_{l\in L\to \infty} q_\rho(x_l,y_l) \\
&= q_\rho(x_*,y_*).
\end{aligned}
\end{equation*}
By the continuity of  $q_{\rho}(x_{l},y_{l})$ and taking the limit of both sides of (\ref{eqn: q_xy-inequality}) as $l\in L\to \infty$, we have 
\begin{equation*} 
  q_{\rho}(x_{*},y_{*}) \le q_{\rho}(x,y_{*}), \quad \forall x\in \mathcal X, \quad \mbox{and} \quad 
  q_{\rho}(x_{*},y_{*}) \le  q_{\rho}(x_{*},y), \quad \forall y\in \mathcal Y.
\end{equation*}
Further, it is clear from  (\ref{eqn: non-increasing_q}) that $\left\{q_{\rho}(x_{l},y_{l})\right\}$ is non-increasing. 

Next, suppose $q_\rho(x_{r}, y_{r}) = q_\rho(x_{r+1}, y_{r+1})$ for some $r \in \mathbb N$. Then by (\ref{eqn: q_xy-inequality}), we have $q_\rho(x_{r+1}, y_{r}) = q_\rho(x_{r}, y_{r})$. Furthermore, since $x_{r+1} \in \mathrm{Argmin}_{x\in \mathcal X} \ q_\rho(x, y_{r})$, we see $q_\rho (x_{r+1}, y_{r}) = \min_{x \in \mathcal X} q_\rho(x, y_{r})$. In view of $q_\rho(x_{r+1}, y_{r}) = q_\rho(x_{r}, y_{r})$ and $x_{r} \in \mathcal X$, we have $q_\rho(x_{r}, y_{r})=\min_{x \in \mathcal X} q_\rho(x, y_{r})$ such that $x_{r} \in \mathrm{Argmin}_{x\in \mathcal X} \ q_\rho(x, y_{r})$. Further, $y_{r} \in \mathcal Y$ satisfies $y_{r} \in \mathrm{Argmin}_{y \in \mathcal Y} \ q_\rho(x_{r}, y)$. This shows that $(x_{r}, y_{r})$ is a saddle point of $(P_{x,y})$.

\end{proof}

\subsection{Proof of Theorem \ref{thm: CCMV-PD convergence}}
\begin{proof}
(i)
{\color{black}
By Lemma \ref{lem: BCD boundedness}, we know that $\{(x^{(j)}, y^{(j)})\}$ is bounded and therefore, has a convergent subsequence. For our purposes, without loss of generality, we suppose that the sequence itself is convergent. Let  $\left(x^*, y^*\right)$ be its accumulation point.
Note that:
$$
q_{\rho^{(j)}}(x^{(j)}, y^{(j)})
\le
q_{\rho^{(j)}}(x_1^{(j)}, y_1^{(j)})
\le
q_{\rho^{(j)}}(x_1^{(j)}, y_0^{(j)})
=
\min_{x\in \mathcal{X}} q_{\rho^{(j)}}(x, y_0^{(j)}) \le  \min_{x\in \mathcal{X}} q_{\rho^{(j+1)}}(x, y_0^{(j)}),
$$
where the first and second inequalities are the results of  $\{q_{\rho^{(j)}}(x_l^{(j)}, y_l^{(j)})\}_{l\in \mathbb N}$ being non-increasing (proved in Theorem 2.1), the equality is just definition, and the last inequality is due to $\rho^{(j)}<\rho^{(j+1)}$. Therefore, because of step (12) of Algorithm 2 and the definition of $q_\rho(x,y)$ given in (3), we see that  
$$
f(x^{(j)})+\rho^{(j)}\|x^{(j)}- y^{(j)}\|_2^2=q_{\rho^{(j)}}(x^{(j)},y^{(j)})
\le \min_{x\in \mathcal{X}} q_{\rho^{(j+1)}}(x, y_0^{(j)}) \le \Upsilon.
$$
Hence, we see  $\{x^{(j}\} \subseteq X_{\Upsilon}$, where $X_{\Upsilon}$ is defined in (10). Next, we can say  and 
$$
\|x^{(j)}- y^{(j)}\|_2^2
\le 
\frac{
\Upsilon-f(x^{(j)})}{\rho^{(j)}}
\le 
\frac{\Upsilon-\min_{x\in X_{\Upsilon}}f(x)}{\rho^{(j)}}.$$
Thus, because $\min_{x\in X_{\Upsilon}}f(x)>-\infty$ when $A\succ 0$, 
by letting $j\to \infty,$ we obtain that   $x^*=y^*$.
}

Next, let $\Ical^{(j)}\subseteq [n]$ be such that $|\Ical^{(j)}|=k$ and $(y_{(\Ical^{(j)})^c})_i=0$ for every $j\in \mathbb N$ and $i\in (\Ical^{(j)})^c$. Then, since $\{\Ical^{(j)}\}$ is a bounded sequence of indices, it has a convergent subsequence, which means that there exists an index subset $\Lcal \subseteq [n]$ with $|\Lcal|=k$ and a subsequence $\{(x^{(j_\ell)}, y^{(j_\ell)}\}$ of the above convergent subsequence such that $\Ical^{(j_\ell)} = \Lcal$ for all large $j_\ell$'s. Therefore, since $x^*=y^*$ and $y^*_{\Lcal^c}=0$, we see $x^*_{\Lcal^c}=0$. Further, $x^*_{\Lcal}\ge 0$ because 
for each $j_\ell$ and $i\in \Ical^{(j_\ell)}$, we know $ (y^{(j_\ell)})_i\ge 0$ such that $ y^*_{\Lcal}\ge 0$ and $x^*_{\Lcal}= y^*_{\Lcal}$.

(ii)  For each $j$, $(x^{(j)}, y^{(j)})$ is a saddle point of  (\ref{pr: pxy}) with $\rho=\rho^{(j)}>0$, so we have
\begin{equation*} \label{eqn: KKT -pxyz-convergence}
\left\{ \begin{split}
& 2Ax^{(j_{\ell})}-\tau \mu +2\rho^{(j_{\ell})}(x^{(j_{\ell})}-y^{(j_{\ell})})+\beta^{(j_\ell)} e=0,
\\ 
&
y^{(j_{\ell})}=H_k(\max(x^{(j_{\ell})},0)= (\max(x^{(j_{\ell})},0))_{\Lcal}
\\  
&
e^T x^{(j_{\ell})}=1, \quad  (y^{(j_{\ell})})_\Lcal\ge 0, \quad \mbox{and} \quad (y^{(j_{\ell})})_{\Lcal^c}=0.
\end{split} \right.
\end{equation*}
By the second equation above, one can see $\left(x^{(j_{\ell})}-y^{(j_{\ell})}\right)_\Lcal=
\left(x^{(j_{\ell})}-\max(x^{(j_{\ell})},0)\right)_\Lcal
= -
( \max(-x^{(j_{\ell})},0))_\Lcal
$
such that $ \max(-x^{(j_{\ell})},0))_\Lcal  \perp (y^{(j_\ell)})_\Lcal$. 
Hence, we get 
\begin{eqnarray}  \label{eqn:kkt-limit}
2Ax^{(j_{\ell})}-\tau\mu -
\underbrace{\begin{bmatrix}  2\rho^{(j_\ell)}( \max(-x^{(j_{\ell})},0))_\Lcal \\ 0\end{bmatrix}}_{:=\lambda^{(j_\ell)}
} 
+
\underbrace{\begin{bmatrix}  0 \\ 2\rho^{(j_\ell)}(  x^{(j_\ell)})_{\Lcal^c} \end{bmatrix}  }_{:=w^{(j_\ell)}}+\beta^{(j_\ell)}e=0,
\end{eqnarray}
where $(w^{(j_\ell)})_\Lcal=0$ for each $j_\ell$ and also, $0\le \lambda^{(j_{\ell})}\perp(y^{(j_\ell)})_\Lcal \ge 0,$ and $e^Tx^{(j_{\ell})}=1.$
We next prove that $\{(\lambda^{(j_\ell)},\beta^{(j_\ell)},w^{(j_\ell)})\}$ is bounded  under Robinson's condition on $x^*$. 
 Suppose not, consider the normalized sequence
\[
( \wt \lambda^{(j_\ell)}, \wt \beta^{(j_\ell)}, \wt w^{(j_\ell)}):= \frac{ (\lambda^{(j_\ell)}, \beta^{(j_\ell)}, w^{(j_\ell)})} { \| (\lambda^{(j_\ell)}, \beta^{(j_\ell)}, w^{(j_\ell)}) \|_2}, \qquad \forall \ {j_\ell}.
\]
Through boundedness of this normalized sequence, it has a convergent subsequence $( \wt \lambda^{(j_\ell)}, \wt \mu^{(j_\ell)}, \wt w^{(j_\ell)})$ whose limit is given by $(\wt \lambda_*, \wt \mu_*, \wt w^*)$ such that $\| (\wt \lambda_*, \wt \mu_*, \wt w^*) \|_2=1$. Thus, by passing the limit and boundedness of $2Ax^{(j_\ell)}-\tau \mu$, we have
\begin{equation*} \label{eqn:limit_condition}
   -\wt\lambda_*  + \wt w^*  + \wt \beta_* e=0,
\end{equation*}
where $\wt \lambda_* \ge 0$ and $\wt w^*_\Lcal=0$. 
By Robinson's conditions at $x^*$, there exist  vectors $d \in \mathbb R^n$ and  $v\in \mathbb R^n$ with $v_i\le 0$ for $i\in \mathcal L^c$ and $\beta \in \mathbb R$ such that $-d - v = -\wt \lambda_*$, $ e^T d = -\wt \beta_*$, and $d_{\Lcal^c} = - \wt w^*_{\Lcal^c}$. Since $d_{\Lcal^c} = - \wt w^*_{\Lcal^c}$ and $\wt w^*_\Lcal=0$, we see that $d^T \wt w^* = -\| \wt w^* \|^2_2$. Therefore,
\[
  0 = -d^T\wt \lambda_* + d^T \wt w^* + \wt \beta_* d^T e=  - (\wt \lambda_*)^2 + \frac{\wt \lambda_* v }{2} - (\wt \beta_*)^2 - \|\wt w^* \|^2_2 = - \|  (\wt \lambda_*, \wt \beta_*, \wt w^*) \|^2_2 +  \frac{\wt \lambda_* v }{2},
\]
which implies that $\|  (\wt \lambda_*, \wt \beta_*, \wt w^*) \|^2_2 =  \frac{\wt \lambda_* v }{2}$.
Since $\wt \lambda_* \ge 0$ and $v \le 0$, we have $\|  (\wt \lambda_*, \wt \beta_*, \wt w^*) \|^2_2=0$, which is a  contradiction. Therefore, the sequence $\big( (\lambda^{(j_\ell)}, \beta^{(j_\ell)}, w^{(j_\ell)}) \big)$ is bounded.
Hence, $\{(\lambda^{(j_\ell)},\beta^{(j_\ell)},w^{(j_\ell)})\}$ is bounded and has a convergent subsequence with the limit $(\lambda,\beta, w)$. Thus, through passing limit in (\ref{eqn:kkt-limit}) and applying the results of part (i),  we have the following:
$$ 2Ax^*-\tau \mu -\lambda+w=0, \quad e^Tx^*=1, \quad (x^*)_{\Lcal^c}=0, \quad w^*_\Lcal =0, \quad  0\le \lambda \perp x^*\ge 0 $$ 
where. This means that $x^*$ satisfies the first-order optimality conditions of (\ref{pr: p-original}) given in (\ref{eqn: KKT _conditions for p}). Next, according to Theorem 2.3 in \autocite{lu2010penalty}, since all the objective function is convex and the constraints in this problem (except the sparsity requirement) are linear, $x^*$ is indeed a local minimizer of (\ref{pr: p-original}).
\end{proof}


\end{document}